\documentclass[10pt, a4paper]{article}
\usepackage{amsfonts, amssymb, amsthm, amsmath}

\newcommand{\R}{\mathbb{R}}
\newcommand{\C}{\mathbb{C}}
\newcommand{\Z}{\mathbb{Z}}
\newcommand{\N}{\mathbb{N}}

\newcommand{\Q}{\mathbb{Q}}

\newcommand{\A}{\mathbb{A}}
 
\newcommand{\ma}{\mathbf}
\newcommand{\ben}{\begin{enumerate}}
\newcommand{\een}{\end{enumerate}}
\newcommand{\bit}{\begin{itemize}}
\newcommand{\eit}{\begin{itemize}}
\newcommand{\beq}{\begin{equation}}
\newcommand{\eeq}{\end{equation}}

\newcommand{\ve}{\varepsilon}
\newcommand{\eps}{\epsilon}
\newcommand{\mcal}{\mathcal}

\newcommand{\lab}{\label}

\newcommand{\la}{\lambda}

\newtheorem{thm}{Theorem}
\newtheorem{lem}{Lemma}

\newtheorem*{hyp}{Hypothesis $[\delta, \theta_\delta]$}
\theoremstyle{definition}
\newtheorem*{ack}{Acknowledgement}

\title{Equal sums of like polynomials}

\author{T.D. Browning\\ 
\small\emph{Mathematical Institute,
24--29 St. Giles', 
Oxford OX1 3LB}\\
\small{browning@maths.ox.ac.uk}}
\date{}

\begin{document}

\maketitle

\begin{abstract}
Let $f \in \Z[x]$ be a  polynomial of degree $d$.   We
establish the paucity of non-trivial positive integer solutions to the equation 
$$
f(x_1)+f(x_2)=f(x_{3})+f(x_{4}),
$$
provided that $d \geq 7$. We also investigate the corresponding
situation for equal sums of three like polynomials. 
\end{abstract}

\section{Introduction}

Let $f \in \Z[x]$ be a polynomial of degree $d\geq 3$, and let $B \geq
1$.  Then for any $s \geq 2$ we define  $M_s(f;B)$ to be the number 
of positive integers $x_1,\ldots,x_{2s} \leq B$ such that 
\beq\lab{Xs}
f(x_1)+\cdots+f(x_s)=f(x_{s+1})+\cdots+f(x_{2s}). 
\eeq
It is conjectured that
$M_s(f;B)$ is dominated by the $s!B^s$ trivial solutions, in which
$x_1,\ldots,x_s$ are a permutation of $x_{s+1},\ldots,x_{2s}$.
We therefore 
define $M_s^{(0)}(f;B)$ to be the number of non-trivial solutions to (\ref{Xs}).
It is in the special case $f_0(x)=x^d$ that this
quantity  has received the most attention.  Thanks
to the work of numerous authors it is well known that 
$$
M_2^{(0)}(f_0;B)=o(B^2)
$$ 
for any $d \geq 3$.  Moreover, recent joint work of the author with
Heath-Brown \cite{bhb} has established the bound 
$$
M_3^{(0)}(f_0;B)=o(B^3)
$$
for $d \geq 33$.

Returning to the more general quantity $M_s^{(0)}(f;B)$, it is only in the case
$s=2$ and $d=3$ that the paucity of non-trivial solutions has been
established.  Indeed, using quite elementary means Wooley \cite{wooley} has shown that
$$
M_2^{(0)}(f;B)=O_{\ve,f}(B^{5/3+\ve}),
$$ 
for any cubic polynomial $f$ and any choice of $\ve>0$.  
The best that is known for polynomials of higher degree is the estimate
$M_2^{(0)}(f;B)=O_{\ve,f}(B^{2+\ve})$, that follows
from the trivial estimate $d(n)=O_\ve(n^\ve)$ for the divisor function.  Our first result
rectifies this situation somewhat for polynomials of sufficiently
large degree.

\begin{thm}\lab{equal1}
Let
$f \in \Z[x]$ be a polynomial of degree $d\geq 4$.  Then we have 
$$
M_2^{(0)}(f;B) \ll_{\ve,f} B^{1+\ve}\left(B^{1/3} +
  B^{2/\sqrt{d}+1/(d-1)} \right),
$$
for any $\ve >0$.  In particular we have
$$
M_2(f;B) = 2B^2(1+o(1))
$$ 
for $d \geq 7$, and 
$$
M_2^{(0)}(f;B)\ll_{\ve,f} B^{4/3+\ve}
$$ 
for $d \geq 42$.  
\end{thm}

Apart from in the special case $f=f_0$ and $s=3$,
nobody has yet been able to establish the
paucity of non-trivial solutions to (\ref{Xs})
for any value of $s\geq 3$.  The best available
estimates arise from applications of Hua's inequality \cite{hua}.  In the
case $s=3$ this provides the estimate
$$
M_3(f;B)=O_{\ve,f}(B^{7/2+\ve}),
$$
provided that $f$ has
degree $d \geq 3$.  
Our second result improves upon this estimate for polynomials of sufficiently
large degree.

\begin{thm}\lab{equal2}
Let $f \in \Z[x]$ be a polynomial of degree $d\geq 4$.  Then we have 
$$
M_3^{(0)}(f;B) \ll_{\ve,f} B^{3+\ve}\left(B^{1/3} +
  B^{2/\sqrt{d}+1/(d-1)} \right),
$$
for any $\ve >0$. In particular we have
$$
M_3(f;B) = o(B^{7/2}) 
$$ 
for $d \geq 20$.
\end{thm}

It would be of considerable interest to increase the range of
$d$ in Theorem \ref{equal2}'s estimate for $M_3(f;B)$, for such progress
would have applications to smooth Weyl differencing as
employed by Vaughan and Wooley \cite{v-w} to obtain further
improvements in Waring's problem.  In order to facilitate future
investigations, it will be convenient to state 
the following hypothesis for given $\delta\in \N$ and $\theta_{\delta} \in \R$.

\begin{hyp}
Let $S \subset \A^3$ be a non-singular affine surface of degree
$\delta\geq 2,$ and
let $S_{\delta} \subseteq S$ denote the subset formed by deleting all of the
curves of degree at most $\delta-2$ from $S$, that are defined over
$\Q$.  Then we have 
$$
\#\{(x_1,x_2,x_3) \in S_\delta\cap\Z^3: \, \max \{|x_1|,|x_2|,|x_3|\}
\leq B\}  
\ll_{\ve,\delta} B^{\theta_{\delta}+\ve}.
$$
\end{hyp}

It is worth underlining that the implied constant in this estimate is
assumed to be independent of the coefficients of the polynomial
defining $S$.
It follows from an easy induction argument that Hypothesis $[\delta, 2]$
always holds, although we can actually do
rather better than this.  Indeed one can apply a result of Heath-Brown
\cite[Theorem 14]{annal} along much the same
lines as in the proof of \cite[Equation (1.15)]{annal}, to deduce that
Hypothesis $[\delta, 2/\sqrt{\delta}+1/(\delta-1)]$ holds.  In order to prove Theorems
\ref{equal1} and \ref{equal2} it will therefore suffice to establish
the following result.

\begin{thm}\lab{equal3}
Assume that Hypothesis $[d,\theta_d]$ holds, and let $s \geq 2$.  Then for any
polynomial $f \in \Z[x]$  of degree $d\geq 4$, we have 
$$
M_s^{(0)}(f;B) \ll_{\ve,f} B^{2s-3+\ve}\left(B^{1/3} +
  B^{\theta_d} \right).
$$
\end{thm}

It should be highlighted that Theorem 
\ref{equal3}  is only
interesting for $s=2$ and $s=3$, since for larger values of $s$ it is
beaten by Hua's inequality. We end this section by discussing the
value of the implied constant in Theorem \ref{equal3}.  As it stands
the constant is clearly allowed to depend upon $f$ in some way.
Suppose once and for all that 
\beq\lab{f}
f(x)=a_0x^d+ a_1x^{d-1}+\cdots+a_d,
\eeq 
with $a_0,\ldots,a_d\in \Z$ and $a_0 > 0$.
Then it is in fact possible to
prove that the implied constant is independent of the coefficients of $f$ whenever $a_1=0$,
and that it depends at most upon $d$ and the choice of $\ve>0$.  We
shall content ourselves with merely indicating at which points of the
argument this sort of finer inequality can be retrieved.

\begin{ack}
This work was undertaken while the second author was supported by 
EPSRC grant number GR/R93155/01.
\end{ack}

\section{Proof of Theorem \ref{equal3}}

Suppose that $f$ is given by (\ref{f}) and that $d \geq 4$.  
We begin the proof of Theorem \ref{equal3} by noticing that
$$
a_0^{d-1}d^df(x)=(a_0dx+a_1)^d+b_2(a_0dx+a_1)^{d-2}+\cdots+b_{d-1}(a_0dx+a_1)+b_d,
$$
for appropriate $b_i \in \Z$ depending upon $a_i$ and $d$.
Upon making the change of variables 
$y_i=a_0d x_i+a_1$  for $1 \leq i \leq 2s$,
we therefore obtain the equation
\beq\lab{Ys}
g(y_1)+\cdots+g(y_s)=g(y_{s+1})+\cdots+g(y_{2s})
\eeq
from (\ref{Xs}), where
\beq\lab{g}
g(y)=y^d+b_2y^{d-2}+\cdots+b_{d-1}y.
\eeq
This transformation has the effect of taking positive integer points not exceeding $B$ to positive integer
points  not exceeding $O_{{f}}(B)$.  We stress that this is the main point
of the argument at which a necessary dependence upon the coefficients of $f$
appears.  Such a transformation is in fact unnecessary if $a_1=0$ in
(\ref{f}).  We have therefore shown that
\beq\lab{ineq1}
M_s^{(0)}(f;B) \leq M_s^{(0)}(g;cB),
\eeq
for some constant  $c>0$ depending only upon $f$.

During the course of our argument, it will be
necessary to handle the contribution from certain ``almost trivial''
solutions to (\ref{Ys}) separately.  Let $S_s(B)$ denote the
contribution to $M_s^{(0)}(g;cB)$ from those 
$y_1,\ldots,y_{2s}$ for which
$$
\{y_1,\ldots,y_s\} \cap \{y_{s+1},\ldots,y_{2s}\} = \emptyset,
$$
and let $T_s(B)$ denote the remaining contribution.  It
follows that 
\beq\lab{ineq2}
M_s^{(0)}(g;cB)=S_s(B)+ T_s(B).
\eeq
Moreover, whenever the vector $(y_1,\ldots,y_{2s})$ is counted by
$T_s(B)$, we must have $y_i=y_j$ for some $1 \leq i \leq s< j \leq
2s$.

In order to estimate $S_s(B)$ and $T_s(B)$ we shall employ the
following result due to Pila \cite[Theorem A]{pila95}.  

\begin{lem}\lab{pila}
Let $C \subset \A^3$
be an absolutely irreducible affine curve of degree $\delta$.  Then we have
$$
\#\{(x_1,x_2,x_3) \in C\cap \Z^3: \, \max \{|x_1|,|x_2|,|x_3|\}
\leq B\} \ll_{\ve,\delta} B^{1/\delta+\ve}.
$$
\end{lem}

As in the
statement of Hypothesis $[\delta,\theta_\delta]$ the implied constant
in Lemma \ref{pila} is understood
to be independent of the coefficients of the polynomials defining $C$.

\subsection{Estimating $S_s(B)$}

In this section we provide an estimate for $S_s(B)$.  This constitutes the
main part of our argument.
The idea will simply be to count points on the affine surfaces obtained by
fixing values of $y_{4}, \ldots, y_{2s}$ in (\ref{Ys}).  
Let 
$$
\eps_i=
\left\{
\begin{array}{ll}
-1, & i \leq s,\\
+1, & i >s,
\end{array}
\right.
$$
and write  $\mcal{N}$ for the set of vectors
$\ma{n}=(n_4, \ldots,n_{2s}) \in (\N\cap[1,cB])^{2s-3}$.
For any $\ma{n} \in \mcal{N}$ we define the surface
$$
\Gamma_{\ma{n}}: \quad g(y_1)+g(y_2)-\eps_3g(y_3)=\sum_{i=4}^{2s}\eps_i g(n_i).
$$
Let $\mcal{N}_1$ be the set of $\ma{n}\in\mcal{N}$ for which 
$\Gamma_\ma{n}$ is singular, and let $\mcal{N}_{2}=\mcal{N}\setminus \mcal{N}_1$.
Clearly $\Gamma_\ma{n}$ is non-singular, and so absolutely
irreducible,  for $\ma{n} \in \mcal{N}_2$.
Our first task is to establish the following result, which ensures
that the same is true for $\ma{n} \in \mcal{N}_1$.

\begin{lem}\lab{N_s}
The surface $\Gamma_\ma{n}$ is absolutely irreducible for any $\ma{n}
\in \mcal{N}_1$, and we have
$\#\mcal{N}_1 =O_{f}(B^{2s-4})$.
\end{lem}
\begin{proof}
Suppose that $(\xi_1,\xi_2,\xi_3)$ is a singular point of the surface
$\Gamma_\ma{n}$, for any $\ma{n} \in \mcal{N}_1$.  Then it follows that
$\frac{dg}{dx}$ vanishes at $\xi_i$ for $1 \leq i \leq 3,$ and that
\beq\lab{xi}
g(\xi_1)+g(\xi_2)-\eps_3g(\xi_3)=\sum_{i=4}^{2s}\eps_i g(n_i).
\eeq
Since $\frac{dg}{dx}$ is a polynomial of degree $d-1$, it follows that there
are at most $(d-1)^3$ possible singular points $(\xi_1,\xi_2,\xi_3)
\in \C^3$
on the surface $\Gamma_\ma{n}$. 
This establishes that $\Gamma_\ma{n}$ is absolutely irreducible.
Indeed, if we had a non-trivial decomposition of the form
$$
\Gamma_\ma{n}=\Gamma_1\cup\Gamma_2\subset \A^3,
$$ 
then $\Gamma_1$ and $\Gamma_2$ would necessarily 
intersect in a variety of dimension  at least $1$.
Since every point of this set 
would produce a singular point in the surface $\Gamma_\ma{n}$,
this would contradict the fact that it has finite singular locus.

It therefore  remains to count the number
of $\ma{n} \in \mcal{N}$ for which (\ref{xi}) holds, for $O_{d}(1)$
values of $(\xi_1,\xi_2,\xi_3)$.  But if we fix a choice of
$(\xi_1,\xi_2,\xi_3)$ and $(n_5,\ldots,n_{2s})$, then there can
clearly only
be $O_{d}(1)$ values of $n_4$ such that (\ref{xi}) holds.  Hence it
follows that there are $O_{f}(B^{2s-4})$ values of $\ma{n} \in \mcal{N}$
such that $\Gamma_{\ma{n}}$ has singularity  $(\xi_1,\xi_2,\xi_3)$.
This suffices to establish the second part of the lemma.
\end{proof}

Let $S_s(B;\ma{n})$ denote the number of positive integers  
$y_1,y_2,y_3 \ll_{f} B$, that lie on the surface $\Gamma_\ma{n}$, with
the constraint that $y_3 \not\in \{y_1,y_2\}$ whenever $s=2$.
Then in order to estimate $S_s(B)$, it will suffice to estimate $S_s(B;\ma{n})$ 
for  each $\ma{n}\in \mcal{N}=\mcal{N}_1 \cup \mcal{N}_2$, since we clearly have
\beq\lab{belfast}
S_s(B) \leq \sum_{\ma{n}\in \mcal{N}_1} S_s(B;\ma{n})
+ \sum_{\ma{n}\in \mcal{N}_2} S_s(B;\ma{n}).
\eeq
In estimating $S_s(B;\ma{n})$ it will prove necessary to pay
special attention to the points lying on curves of low degree
contained in $\Gamma_\ma{n}$.

\begin{lem}\lab{low}
For any $\ma{n}\in \mcal{N}$, there is no contribution to
$S_s(B;\ma{n})$ from any lines or conics contained in
$\Gamma_\ma{n}$ that are defined over $\Q$.
\end{lem}

\begin{proof}
We begin by considering the
possibility that $\Gamma_\ma{n}$ contains a line defined over $\Q$, and
we write $c_\ma{n}=\sum_{i=4}^{2s}\eps_i g(n_i)$
for convenience.
Thus there exist $\la_i, \mu_i \in \Q$  such that
the polynomial
$$
g(\la_1 t +\mu_1)+g(\la_2 t +\mu_2)-\eps_3g(\la_3 t +\mu_3)-c_\ma{n}
$$
vanishes identically.  We may clearly assume that at most one $\la_i$ is
zero, since otherwise it is easy to see that $\la_1=\la_2=\la_3=0$.
Suppose first that $\la_1 \neq 0$.  Then after a possible change of
variables we may assume that $\la_1=1$ and $\mu_1=0$.  Upon recalling
the shape (\ref{g}) that $g$ takes, it therefore follows that
$$
1+\la_2^d=\eps_3 \la_3^d. 
$$
Wiles' proof \cite[Theorem 0.5]{flt} of Fermat's Last Theorem shows that
$\la_2\la_3=0$. If $\la_3=0$, then $d$ must be odd and $\la_2=-1$.  We
must now consider the possibility that we have an identity of the
shape $g(t)+g(-t+\mu_2)=k_\ma{n}$, for some constant $k_\ma{n}$.  Upon examining the
coefficient of $t^{d-1}$, we conclude from (\ref{g}) that $\mu_2=0$.
In terms of the original coordinates we have shown that  this case 
produces the affine line $y_1=-y_2, y_3=\mu_3$, provided that $d$ is odd.
Although this line may be contained in $\Gamma_\ma{n}$ for certain
choices of $g$, such solutions actually contribute nothing to
$S_s(B;\ma{n})$ since we are only interested in
positive integer points on $\Gamma_\ma{n}$.
Next we suppose that $\la_2=0$ and 
$\la_3^d=\eps_3$.  If $\eps_3=-1$ then the previous argument can be
repeated to yield the line $y_1=-y_3, y_2=\mu_2$, which shows
that this case also contributes nothing to $S_s(B;\ma{n})$.  If $\eps_3=1$
however, so that $s=2$, then either $\la_3=-1$ and $d$ is even, or
else $\la_3=1.$ In the former case we obtain the line $y_1=-y_3,
y_2=\mu_2$, and in the latter case we obtain the line $y_1=y_3, y_2=\mu_2$.
Neither of these cases contribute anything to $S_2(B;\ma{n})$, since
we must have $y_3 \not\in \{y_1,y_2\}$ whenever $s=2$.
Upon treating the case corresponding to $\la_1=0$ in a similar
fashion, one is led to the pair of lines $y_1=\mu_1, y_2=\pm
y_3$.  Neither of these contribute anything to
$S_s(B;\ma{n})$.

Next we suppose that $\Gamma_\ma{n}$ contains a 
conic defined over $\Q$.   Thus 
there exist $\kappa_i, \la_i, \mu_i \in \Q$ such that the polynomial
$$
g(\kappa_1t^2+\la_1 t +\mu_1)+g(\kappa_2t^2+\la_2 t +\mu_2)-\eps_3g(\kappa_3t^2+\la_3 t +\mu_3)-c_\ma{n}
$$
vanishes identically.  
Suppose that $\kappa_1 \neq 0$, say, and let 
$$
h(y)=\kappa_1^d g(y/\kappa_1)=y^d+b_2\kappa_1^2y^{d-2}+\cdots+b_{d-1}\kappa_1^{d-1}y.
$$
Then after a possible relabling of variables it suffices to consider
the vanishing of the polynomial
$$
h(t^2+\mu_1)+h(\kappa_2t^2+\la_2 t
+\mu_2)-\eps_3h(\kappa_3t^2+\la_3 t +\mu_3)-\kappa_1^dc_\ma{n}.
$$
Upon examining the coefficient of the leading monomial $t^{2d}$, we deduce that
$$
1+\kappa_2^d=\eps_3 \kappa_3^d,
$$
and so $\kappa_2\kappa_3=0$.  If   $\kappa_3=0$ then $d$ must be odd and $\kappa_2=-1$. 
Using the fact that the coefficient of $t^{2d-1}$ must also vanish, we
further deduce that 
$\la_2=0$.    Similarly, the coefficient of $t^{2d-2}$ is equal to 
$d(\mu_1+\mu_2)$  since $d$ is odd, 
from which it follows that $\mu_2=-\mu_1$.
Again appealing to the fact that $d$ is odd, we finally deduce that
$\la_3=0$ by considering the coefficient of $t^{d}$.
In terms of the original coordinates we therefore have
$$
y_1=t^2+\mu_1=-y_2, \qquad y_3=\mu_3,
$$
and it is clear that such solutions contribute nothing to $S_s(B;\ma{n})$.
Alternatively, if 
$\kappa_2=0$ then $\kappa_3^d=\eps_3$.  Arguing as above it suffices
to assume that $\eps_3=1$, and so that $s=2$.  But then the same
analysis ultimately leads to solutions of the form
$$
y_1=t^2+\mu_1=y_3, \qquad y_2=\mu_2,
$$
which are not permissible.  The case $\kappa_1=0, \kappa_2 \neq 0$ is
handled similarly. This completes the proof of Lemma \ref{low}.
\end{proof}

The remainder of this section is taken up with establishing the
following result.

\begin{lem}\lab{surface}
Assume that Hypothesis $[d,\theta_d]$ holds.  Then we have
$$
S_s(B;\ma{n}) \ll_{\ve,f} 
\left\{ 
\begin{array}{ll}
B^{2/\sqrt{d}+1/3+\ve}, & \ma{n} \in \mcal{N}_{1},\\
B^{1/3+\ve}+B^{\theta_d+\ve}, & \ma{n} \in \mcal{N}_{2}.\\
\end{array}
\right.
$$
\end{lem}

\begin{proof}
The proof of Lemma \ref{surface} will hinge upon work of
Heath-Brown \cite[Theorem 14]{annal}.  Let $\|\Gamma_\ma{n}\|$ denote the maximum
modulus of the coefficients of the polynomial defining
$\Gamma_\ma{n}$, so that in particular 
$\log \|\Gamma_\ma{n}\| =O_{f} (\log B)$ for any $\ma{n} \in
\mcal{N}_1\cup \mcal{N}_2$.
Now let $(y_1,y_2,y_3)$ be any point counted by
$S_s(B;\ma{n})$, and note that $\Gamma_\ma{n}$ is absolutely irreducible for $\ma{n} \in
\mcal{N}_1\cup \mcal{N}_2$ by Lemma \ref{N_s}.
Therefore an application of 
\cite[Theorem 14]{annal} implies that
$(y_1,y_2,y_3)$ lies on one of at most 
$$
\ll_{\ve,f}B^{2/\sqrt{d}+\ve}\|\Gamma_\ma{n}\|^\ve
\ll_{\ve,f}B^{2/\sqrt{d}+\ve}
$$
proper
subvarieties of $\Gamma_\ma{n}$, each of degree $O_{\ve,d}(1)$.
We remark that it is actually possible to make the assumption 
$\log \|\Gamma_\ma{n}\| =O_{\ve,d}(\log B)$ at this stage, by
employing the argument of \cite[Theorem $4$]{annal}.  When $a_1=0$ in
(\ref{f}) this would lead to the uniformity result mentioned after the
statement of Theorem \ref{equal3}.

It remains to estimate the number of points of bounded height lying on
$O_{\ve,f}(B^{2/\sqrt{d}+\ve})$ absolutely irreducible curves of degree
$O_{\ve,d}(1)$ that are contained in 
$\Gamma_\ma{n}$.   For this we use Lemma \ref{pila}.
Suppose first that $\ma{n} \in \mcal{N}_1$.  Then Lemma \ref{low}
implies that we may ignore points lying on any curves of degree at
most $2$ that are defined over $\Q$. 
Any curve of degree at most $2$ that is not defined over $\Q$ clearly
contains only $O(1)$ points.
Hence it follows from Lemma \ref{pila} that 
$$
S_s(B;\ma{n}) \ll_{{\ve,f}} B^{2/\sqrt{d}+1/3+\ve},
$$
whenever  $\ma{n} \in \mcal{N}_1$.

Suppose now that $\ma{n} \in \mcal{N}_2$. Then on the assumption that
Hypothesis $[d,\theta_d]$ holds, we obtain the overall contribution $O_{\ve,f}(B^{\theta_{d}+\ve})$ from
points not contained on any curve of degree at most $d-2$ that is
defined over $\Q$.
It remains to consider the  contribution to $S_s(B;\ma{n})$ from the
curves of degree at most $d-2$,  that are
defined over $\Q$ and are contained in $\Gamma_\ma{n}$.
Since $\Gamma_\ma{n}$ is non-singular we may apply a result of 
Colliot-Th\'el\`ene \cite[Appendix]{annal}.  
We conclude that $\Gamma_\ma{n}$ contains  $O_{d}(1)$ curves of degree $\leq
d-2$. Lemma \ref{low} implies that we may ignore points lying on those curves of degree 
at most $2$.  Hence Lemma \ref{pila} yields the overall contribution
$O_{\ve,f}(B^{1/3+\ve})$  from the curves of degree at
most $d-2$ contained in $\Gamma_\ma{n}$.  
This completes the proof of Lemma \ref{surface}.
\end{proof}

Recall the estimate in Lemma \ref{N_s} for $\#\mcal{N}_1$, and note
that $\#\mcal{N}_2=O_{f}(B^{2s-3})$.  Then we may combine 
Lemma \ref{surface} and  (\ref{belfast}) to deduce that
\beq\lab{S-final}
S_s(B) \ll_{\ve,f} B^{2s-3+\ve}\left(B^{1/3} +
  B^{\theta_d} \right),
\eeq
since $d \geq 4$.

\subsection{Estimating $T_s(B)$}

In this section we shall study the quantity $T_s(B)$.  Under the
assumption that Hypothesis $[d,\theta_d]$ holds, our aim is to establish the inequality
\beq\lab{a-triv}
T_s(B)\ll_{\ve,f} B^{2s-3+\ve}\left(B^{1/3} +
  B^{\theta_d} \right),
\eeq
for any $s\geq 2$.  We shall argue by induction on $s$.  

In order to handle the base
case $s=2$, it will suffice to estimate the contribution to $T_2(B)$ from
those $y_1,y_2,y_3,y_4$ for which $y_1=y_3$, say.  There are
then $O_{f}(B)$ choices for $y_1,y_3$, and it remains to count the
number of positive integers $y_2,y_4 \ll_{f} B$ such that $y_2 \neq
y_4$ and 
\beq\lab{curve}
g(y_2)=g(y_4).
\eeq
This equation defines a curve of degree $d$ in $\A^2$. We claim that
those points lying on curves of degree at most $2$, that form
components of (\ref{curve}) and are defined over $\Q$, contribute nothing to
$T_2(B)$.  This is established along exactly the same lines as the proof
of Lemma \ref{low}, and so we will be brief.
Suppose first that there exist $\la,\mu,\la', \mu' \in \Q$  such that
$g(\la t +\mu)=g(\la' t +\mu')$
vanishes identically.  After a possible change of variables we may 
assume without loss of generality that $\la'=1$ and $\mu'=0$.
By equating coefficients it therefore follows from (\ref{g}) that 
$\la=\pm 1$ and $\mu=0$.   Neither of
these cases contribute anything to $T_2(B)$.  The case in which
(\ref{curve}) contains a conic defined over $\Q$ is despatched in
precisely the same way.
Thus it remains to
estimate the contribution from the remaining absolutely irreducuble components
of (\ref{curve}).  An application of Lemma \ref{pila}
therefore yields the overall
contribution $O_{\ve,f}(B^{1/3+\ve})$ to $T_2(B)$.  This 
establishes that $T_2(B)=O_{\ve,f}(B^{4/3+\ve})$, which is
satisfactory for (\ref{a-triv}).

Suppose now that $s>2$.  But then it is trivial to see that we have
$$
T_s(B) \ll \sum_{y\ll_{f} B} \left(S_{s-1}(B)+T_{s-1}(B)\right).
$$
Applying the induction hypothesis, in conjunction with (\ref{S-final}), therefore yields 
$$
T_s(B) \ll_{\ve,f}  \sum_{y\ll_{f} B} B^{2s-5+\ve}\left(B^{1/3} +
  B^{\theta_d} \right) \ll_{\ve,f} B^{2s-3+\ve}\left(B^{1/3} + 
  B^{\theta_d} \right).
$$
This completes the proof of (\ref{a-triv}).

\subsection{Completion of the proof}

Assume that Hypothesis $[d,\theta_d]$ holds.
Then it remains to combine (\ref{S-final}) and (\ref{a-triv}) in (\ref{ineq1})
and (\ref{ineq2}), to conclude that 
$$
M_s^{(0)}(f;B) \leq S_s(B)+ T_s(B) \ll_{\ve,f} B^{2s-3+\ve}\left(B^{1/3} + B^{\theta_d} \right).
$$
This completes the proof of Theorem \ref{equal3}.

\end{document}